\documentclass[11pt, a4paper]{article}

\usepackage{amsthm}
\usepackage{amsmath}
\usepackage{amssymb}
\usepackage{amscd}                 
\usepackage[all]{xy}               
\usepackage{stmaryrd}              
\usepackage{a4wide}
\usepackage{url}
\theoremstyle{plain}

\newtheorem{thm}{Theorem}[section]
\newtheorem{cor}[thm]{Corollary}

\theoremstyle{definition}

\newtheorem{rem}[thm]{Remark}

\newcommand\GSp{\mathrm{GSp}}
\newcommand\PGSp{\mathrm{PGSp}}
\newcommand\Sp{\mathrm{Sp}}
\newcommand\PSp{\mathrm{PSp}}
\newcommand\GL{\mathrm{GL}}
\newcommand\PSL{\mathrm{PSL}}
\newcommand\PGL{\mathrm{PGL}}
\newcommand\PSU{\mathrm{PSU}}
\newcommand\SU{\mathrm{SU}}
\newcommand\SL{\mathrm{SL}}
\newcommand\Gal{\mathrm{Gal}}
\newcommand\Q{\mathbb{Q}}
\newcommand{\id}{\mathrm{id}}

\title{Galois representations and Galois groups over $\mathbb{Q}$}

\author{Sara Arias-de-Reyna, C\'ecile Armana, Valentijn Karemaker,\\ Marusia Rebolledo, Lara Thomas and N\'uria Vila}
\date{}

\begin{document}

\maketitle

\begin{abstract}
 In this paper we generalize results of P.~Le Duff to genus $n$ hyperelliptic curves. More precisely, let
 $C/\mathbb{Q}$ be a hyperelliptic genus $n$ curve, let  $J(C)$ be the associated Jacobian variety and let
 $\bar{\rho}_{\ell}:G_{\mathbb{Q}}\rightarrow \GSp(J(C)[\ell])$ be the Galois representation attached to the $\ell$-torsion
 of $J(C)$.
 Assume that there exists a prime $p$ such that $J(C)$  has semistable 
 reduction with toric dimension 1 at $p$. We provide an algorithm to compute a list
 of primes $\ell$  (if they exist) such that $\bar{\rho}_{\ell}$ is surjective. In particular we realize $\mathrm{GSp}_6(\mathbb{F}_{\ell})$ as a Galois group over $\mathbb{Q}$ for 
 all primes $\ell\in [11, 500000]$.
\end{abstract}

\section*{Introduction}

In this paper we present the work carried out at the conference
\emph{Women in numbers - Europe}, (October 2013), by the working group 
\emph{Galois representations and Galois groups over $\mathbb{Q}$}.
Our aim was to study the image of Galois representations attached to
the Jacobian varieties of genus $n$ curves, motivated by the applications to the 
inverse Galois problem over $\mathbb{Q}$. In the case of genus 2, 
there are several results in this direction (e.g. \cite{LeDuff}, \cite{Di02}), and
we wanted to explore  the scope of these results. 

Our result is a generalization of P.~Le Duff's work to the genus $n$ setting, 
which allows us to produce realizations of groups  $\mathrm{GSp}_6(\mathbb{F}_{\ell})$ as 
Galois groups over $\mathbb{Q}$, for infinite families of primes $\ell$ (with positive
Dirichlet density). These realizations are obtained through the Galois representations $\bar{\rho}_{\ell}$
attached to the $\ell$-torsion points of the Jacobian of a genus $3$ curve.

The first section of this paper contains a historical introduction to the inverse Galois problem and some results obtained in this direction by means
of Galois representations associated to geometric objects. Section 2 presents
some theoretic tools, which we collect to prove a result, valid for a class of abelian varieties $A$ of dimension $n$, that yields primes $\ell$ for which we can ensure surjectivity
of the Galois representation attached to the $\ell$-torsion of $A$ (see
Theorem \ref{thm:explicitsurjectivity}). In Section 3, we focus on hyperelliptic curves and explain the computations that allow us to realize $\mathrm{GSp}_6(\mathbb{F}_{\ell})$ as a Galois group over $\mathbb{Q}$
for all primes $\ell\in [11, 500000]$.

\bigskip

\noindent\textbf{Acknowledgements}

\noindent The authors would like to thank Marie-Jos\'e Bertin, Alina Bucur, Brooke Feigon and Leila Schneps for organizing the WIN-Europe conference which initiated this collaboration. Moreover, we are grateful to the Centre International de Rencontres Math\'ematiques, the Institut de Math\'ematiques de Jussieu and the Institut Henri Poincar\'e
for their hospitality during several short visits.
The authors are indebted to Irene Bouw, Jean-Baptiste Gramain, Kristin Lauter, Elisa Lorenzo, Melanie Matchett Wood, Frans Oort and Christophe Ritzenthaler for several insightful discussions. We also want to thank the anonymous referee for her/his suggestions that helped us to improve this paper.

S.~Arias-de-Reyna and N.~Vila are partially supported by the project MTM2012-33830 of the Ministerio de Econom\'ia y Competitividad of Spain,  C.~Armana  by a BQR 2013 Grant from Universit\'e de Franche-Comt\'e and  M.~Rebolledo  by the ANR Project R\'egulateurs  ANR-12-BS01-0002. L.~Thomas thanks the Laboratoire de Math\'ematiques de Besan\c con for its support.


\section{ Images of Galois representations and the inverse Galois problem }\label{sec:1}

One of the main objectives in algebraic number theory is to understand  the absolute Galois group of the rational field,
$G_{\mathbb{Q}}=\Gal (\overline{ \mathbb{Q}}/ \mathbb{Q} ).$
We believe that we would get all arithmetic information if we knew the structure of $G_{\mathbb{Q}}$.
 This is a huge group, but it is compact with respect to the profinite topology.
Two problems arise in a natural way: on the one hand, the identification  of the  finite quotients  of $G_{\mathbb{Q}}$, and on the other hand, the study of $G_{\mathbb{Q}}$ via its Galois representations.

The inverse Galois problem asks whether, for a given finite group
$G$, there exists a Galois extension $L/ \Q$ with Galois group isomorphic to
$G$. In other words, whether a finite group  $G$ occurs as a quotient of $G_{\mathbb{Q}}$.
As is well known, this is an open problem. The origin of this question can be traced back to Hilbert. In 1892, he proved that the symmetric group $S_n$ and the alternating group $A_n$ are Galois groups over $\Q$, for all $n$.
We also have an affirmative answer to the inverse Galois problem for some other families of finite groups. For instance,  all finite solvable groups and  all  sporadic simple groups, except  the Mathieu group $M_{23}$,   are known to be Galois groups over $\Q$.

A Galois representation is a continuous homomorphism
\[  \rho: G_{\mathbb{Q}} \rightarrow \GL_n(R), \]
where $R$ is a topological ring. Examples for $R$ are $\mathbb{C}$,  $\mathbb{Z}/n\mathbb{Z}$ or $\mathbb{F}_q$ with the discrete topology, and
$\mathbb{Q}_\ell$ with the $\ell$-adic topology. Conjectures by
Artin, Serre, Fontaine-Mazur and Langlands, which have experienced significant progress in recent years,  are connected with these Galois representations.

Since $G_{\mathbb{Q}}$ is compact, the image of $\rho$ is finite when the topology of $R$ is discrete.
As a consequence, images of  Galois representations  yield Galois realizations over $\mathbb{Q}$ of finite linear groups
$$\Gal (\overline{ \mathbb{Q}}^{\ker \rho}/\mathbb{Q})\simeq \rho(G_{\mathbb{Q}})\subseteq \GL_n(R).$$

This  gives us an interesting connection between these two questions and  provides us with a strategy to address the inverse Galois problem.

 Let us assume that  $\rho$ is an $\ell$-adic Galois representation associated to some arithmetic-geometric object. In this case, we have additional information   on the ramification behavior, like the characteristic polynomial of the image of the Frobenius elements at unramified primes or the description of the image of the inertia group at the prime $\ell$. This gives us some control on the image of mod $\ell$ Galois representations in some cases and we can obtain, along the way, families of linear groups over finite fields as Galois groups over $\Q$.

More precisely, let $X/\mathbb{Q}$ be a smooth projective  variety and let
$$\rho_\ell: G_\mathbb{Q} \rightarrow
\GL(H^k_{\textrm{\'et}} (X_{\overline{ \mathbb{Q}}}, \mathbb{Q}_\ell)),$$
be  the  $\ell$-adic Galois representation on the $k$-th  \'{e}tale cohomology.
We know  that:
\begin{itemize}
\item $\rho_\ell$ is unramified  away from
$\ell$ and the primes of bad reduction for $X$,
\item if $p$ is a prime of good reduction and  $p\neq \ell$, the characteristic polynomial of  $\rho_\ell(\mathrm{Frob}_p)$ has coefficients
 in $\mathbb{Z}$, is independent of $\ell$ and its
roots have absolute value $p^{k/2}$.\end{itemize}

Let us consider an attached residual Galois representation
$$\overline{\sigma}_\lambda: G_\mathbb{Q} \rightarrow
\GL_n(\mathbb{F}_{\ell^r}),$$
where $\lambda$ is a prime in a suitable number field, dividing $\ell$ and $r\geq 1$ an integer.
To determine the image of $\overline{\sigma}_\lambda$, we usually need to know the classification of maximal subgroups of $ \GL_n(\mathbb{F}_{\ell^r})$, as well as  a description of the image of the inertia group at $\ell$ and the computation of the characteristic polynomial of  $\overline{\sigma}_\lambda(\mathrm{Frob}_p)$, for some prime of good reduction $p\neq \ell$.

\medskip
Let us summarize the known cases of realizations of finite linear groups as Galois groups over $\Q$, obtained via Galois representations.

In the case of  $2$-dimensional  Galois representations  attached to an elliptic curve $E$ defined over $\Q$ without complex multiplication, we know, by a celebrated result of Serre \cite{Proprietes}, that the associated residual Galois representation  is surjective, for all but finitely many primes.
Moreover, it can be shown that if we take, for example,  the elliptic curve $E$ defined by the Weierstrass equation $ Y^2+Y=X^3-X$, then the attached residual Galois representation  is surjective, for all primes $\ell$. Thus we obtain that the group $\GL_2(\mathbb{F}_{\ell})$ occurs as a Galois group over $\Q$, for all primes $\ell$. Actually  we have additional information in this case: the Galois extension  $\Q(E[\ell])/\Q$ is a Galois realization of  $\GL_2(\mathbb{F}_{\ell})$, and it is unramified away from  $37$ and $\ell$, since $E$ has conductor $37$.

The image of $2$-dimensional  Galois representations,  attached to classical modular forms without complex multiplication, has been studied by Ribet  \cite{R}. The image of the residual Galois representations attached to a normalized cuspidal Hecke eigenform without complex multiplication is as large as possible, for all but finitely many primes $\lambda$. This gives us that the groups $\PSL_2(\mathbb{F}_{\ell^r})$  or
 $\PGL_2(\mathbb{F}_{\ell^r})$ can occur as Galois groups over $\Q$.  Moreover, we have effective control of primes with large image for the mod  $\ell$ Galois representation attached to specific modular forms. This gives us Galois realizations over $\mathbb{Q}$ of the groups
$\PSL_2(\mathbb{F}_{\ell^r})$, $r$ even, and
 $\PGL_2(\mathbb{F}_{\ell^r})$, $r$ odd; $1 \leq r\leq 10 $,
 for explicit infinite families of primes $\ell$, given by congruence conditions on $\ell$ (cf. \cite{R-V}, \cite{D-V00}).

Recently, it has been proven that the groups $\PSL_2(\mathbb{F}_{\ell})$ are Galois groups over $\mathbb{Q}$ for all $\ell>3$, by considering the Galois representations attached to an explicit elliptic surface (see \cite{Z}).

Results on generically large image of compatible systems of
3-dimensional Galois representations associated to  some smooth projective surfaces and to some cohomological modular forms  are obtained in \cite{D-V04}.
The  effective control of primes with large image for the residual 3-dimensional Galois representations attached to  some explicit  examples gives us
 that the groups
$\PSL_3(\mathbb{F}_{\ell})$,
$\PSU_3(\mathbb{F}_\ell)$,
$\SL_3(\mathbb{F}_\ell) $,
 $\SU_3(\mathbb{F}_\ell) $
are Galois groups over  $ \mathbb{Q}$, for explicit infinite families of primes $\ell$ (cf. \cite{D-V04}).

In the case of  $4$-dimensional Galois representations, we have  results on  large image for compatible systems of Galois representations  attached to abelian surfaces $A$ defined over $\mathbb{Q}$ such that $\mathrm{End}_{\overline{\mathbb{Q} }}(A)=\mathbb{Z}$, to Siegel modular forms of genus two and to some pure motives (cf. \cite{LeDuff}, \cite{DKR}, \cite{D02}, \cite{D-V11}). The  effective control of primes with large image in some explicit cases gives us that the groups $\PGSp_4(\mathbb{F}_{\ell})$, for all $\ell >3$; and the groups  $\PGSp_4(\mathbb{F}_{\ell^3})$, $\PSp_4(\mathbb{F}_{\ell^2})$, $\PSL_4(\mathbb{F}_{\ell})$ and $\PSU_4(\mathbb{F}_{\ell})$,
for explicit infinite families of  primes $\ell$, are Galois groups over $\Q$ (cf. \cite{ArVi}, \cite{DKR}, \cite{D02}, \cite{D-V08}).\\

In the next section  we consider the image of residual  Galois representations attached to principally polarized abelian varieties of dimension $n$, which provides Galois realizations over $\Q$ of the general symplectic  group $\GSp_{2n}(\mathbb{F}_\ell)$, for almost all $\ell$. \\

Finally, we remark that, using these methods, we can expect to obtain realizations of the groups
$\PSL_2(\mathbb{F}_{\ell^r})$,  $\PGL_2(\mathbb{F}_{\ell^r})$, $\PGSp_{2n}(\mathbb{F}_{\ell^r})$ and  $\PSp_{2n}(\mathbb{F}_{\ell^r})$  as  Galois groups over $\mathbb{Q}$. In fact, by considering compatible systems of Galois representations attached to certain automorphic forms,  we know (cf. \cite{W}, \cite{D-W}, \cite{KLS}, \cite{ArDiShWi}) that these groups are Galois groups over $\mathbb{Q}$, for  infinitely many integers $r$ and infinitely many primes $\ell$. More precisely, we have:

\begin{itemize}
\item ``Vertical direction": For every fixed prime $\ell$, there are infinitely many positive integers~$r$, such that
$\PSL_2(\mathbb{F}_{\ell^r})$ can be realized as a Galois group over $\mathbb{Q}$. Moreover, for each $n \geq 2$, there are infinitely many positive integers~$r$, such that either $\PGSp_{2n}(\mathbb{F}_{\ell^r})$ or $\PSp_{2n}(\mathbb{F}_{\ell^r})$ are Galois groups over $\mathbb{Q}$ (cf.  \cite{W}, \cite{KLS}).
\item ``Horizontal direction":  For every fixed $r$,  there is a positive density set of primes~$\ell$, such that
$\PSL_2(\mathbb{F}_{\ell^r})$ can be realized as a Galois group over $\mathbb{Q}$.  Moreover, for each $n \geq 2$, there is a set of primes~$\ell$ of positive density for which either $\PGSp_{2n}(\mathbb{F}_{\ell^r})$ or $\PSp_{2n}(\mathbb{F}_{\ell^r})$ are Galois groups over $\mathbb{Q}$ (cf. \cite{D-W}, \cite{ArDiShWi}). \end{itemize}


\section{Galois representations attached to abelian varieties}\label{sec:2}

\subsection{The image of the $\ell$-torsion Galois representation}

Let $A$ be an abelian variety of dimension $n$ defined over $\mathbb{Q}$. 
The set of $\overline{\mathbb{Q}}$-points of $A$ admits a group structure. 
Let $\ell$ be a prime number. Then the subgroup of the $\overline{\mathbb{Q}}$-points of $A$ consisting of all
$\ell$-torsion points, 
which is denoted by $A[\ell]$, is  isomorphic to $(\mathbb{Z}/\ell\mathbb{Z})^{2n}$ 
and it is endowed with a natural action of $G_{\mathbb{Q}}$. Therefore, it gives 
rise to a (continuous) Galois representation 
\begin{equation*}
 \overline{\rho}_{A, \ell}:G_{\mathbb{Q}}\rightarrow \GL(A[\ell])\simeq \GL_{2n}(\mathbb{F}_{\ell}).
\end{equation*}
As explained in Section \ref{sec:1}, we obtain a realization of the image of 
$\overline{\rho}_{A, \ell}$ as a Galois group over~$\mathbb{Q}$.

In this section, we will consider principally polarized abelian varieties, 
i.e. we will consider pairs $(A, \lambda)$, where $A$ is an abelian variety 
(defined over $\mathbb{Q}$) and  $\lambda:A\rightarrow A^{\vee}$ is an isogeny 
of degree 1 (that is, an isomorphism between $A$ and the dual abelian variety 
$A^{\vee}$), induced from an ample divisor on $A$. Not every abelian variety 
$A$ admits a principal polarization $\lambda$ and, when it does, it causes 
certain restrictions on the image of $\overline{\rho}_{A, \ell}$. 

Let $V$ be a vector space of dimension $2n$, which is defined over $\mathbb{F}_{\ell}$ 
and endowed with a symplectic (i.e.~skew-symmetric, nondegenerate) pairing 
$\langle \cdot, \cdot \rangle:V\times V\rightarrow \mathbb{F}_{\ell}$. We consider 
the \emph{symplectic group} 
\begin{equation*}\Sp(V, \langle \cdot, \cdot \rangle):= \{M\in \GL(V): \forall v_1, v_2\in V, \langle Mv_1, Mv_2\rangle=\langle v_1, v_2\rangle\}\end{equation*} 
and the \emph{general symplectic group} 
\begin{equation*}\GSp(V, \langle \cdot, \cdot \rangle):= \{M\in \GL(V): \exists m\in \mathbb{F}_{\ell}^{\times}\text{ such that }\forall v_1, v_2\in V, \langle Mv_1, Mv_2\rangle=m \langle v_1, v_2\rangle\}.\end{equation*} 
When $A$ is a principally polarized abelian variety, the image of $\overline{\rho}_{A, \ell}$ lies 
inside the general symplectic group of $A[\ell]$ with respect to a certain symplectic pairing.
More precisely, denote by $\mu_{\ell}(\overline{\mathbb{Q}})$ the group of $\ell$-th roots 
of unity inside a fixed algebraic closure $\overline{\mathbb{Q}}$ of $\mathbb{Q}$. Recall that the Weil 
pairing $e_{\ell}$ is a perfect pairing
\begin{equation*}
 e_{\ell}:A[\ell]\times A^{\vee}[\ell]\rightarrow \mu_{\ell}(\overline{\mathbb{Q}}).
\end{equation*}
If $(A, \lambda)$ is a principally polarized abelian variety, we can consider the pairing 
\begin{equation*}\begin{aligned}
 e_{\ell, \lambda}:A[\ell]\times A[\ell]&\rightarrow \mu_{\ell}(\overline{\mathbb{Q}})\\
 (P, Q)& \mapsto e_{\ell}(P, \lambda(Q))\end{aligned}
\end{equation*}
which is a non-degenerate skew-symmetric pairing (i.e.~a symplectic pairing), compatible 
with the action of $G_{\mathbb{Q}}$. This last condition means that, for any $\sigma\in G_{\mathbb{Q}}$, 
\begin{equation*}
(e_{\ell, \lambda}(P, Q))^{\sigma}=e_{\ell, \lambda}(P^{\sigma}, Q^{\sigma}).
\end{equation*}

Note that $G_{\mathbb{Q}}$ acts on $\mu_{\ell}(\overline{\mathbb{Q}})$ via the mod $\ell$ 
cyclotomic character $\chi_{\ell}$, so that 
$(e_{\ell, \lambda}(P, Q))^{\sigma}=(e_{\ell, \lambda}(P, Q))^{\chi_{\ell}(\sigma)}$. 
If we fix a primitive $\ell$-th root of unity $\zeta_{\ell}$, we may write the pairing 
$e_{\ell, \lambda}(\cdot, \cdot)$ additively, i.e. we define 
\begin{equation*}\langle \cdot, \cdot \rangle:A[\ell]\times A[\ell]  \rightarrow \mathbb{F}_{\ell} \end{equation*} 
as $\langle P, Q\rangle:= a\text{ such that }\zeta^a=e_{\ell, \lambda}(P, Q)$. 

In other words, we have a symplectic pairing on the $\mathbb{F}_{\ell}$-vector 
space $A[\ell]$ such that, for all $\sigma\in G_{\mathbb{Q}}$, 
the linear map $\overline{\rho}(\sigma):A[\ell]\rightarrow A[\ell]$ satisfies that 
there exists a scalar, namely $\chi_{\ell}(\sigma)$, such that 
\begin{equation}\label{eq:multiplier} \langle\overline{\rho}(\sigma)(P), \overline{\rho}(\sigma)(Q)\rangle = \chi_{\ell}(\sigma)\langle P, Q\rangle.\end{equation}
That is to say, the image of the representation $\overline{\rho}_{A, \ell}$ 
is contained in
the general symplectic group $\GSp(A[\ell], \langle \cdot, \cdot \rangle)\simeq \GSp_{2n}(\mathbb{F}_{\ell})$. 
Therefore, below we will consider $\overline{\rho}_{A,\ell}$ as a map into $\GSp(A[\ell], \langle \cdot, \cdot \rangle)\simeq \GSp_{2n}(\mathbb{F}_{\ell})$ and we will say that it is surjective if $\mathrm{Im}\overline{\rho}_{A, \ell}=\GSp(A[\ell])\simeq \GSp_{2n}(\mathbb{F}_{\ell})$. 

The determination of the images of the Galois representations $\overline{\rho}_{A, \ell}$ 
attached to the $\ell$-torsion of abelian varieties is a topic that has received a lot of 
attention. A remarkable result  by Serre quoted  in \cite[n. 136, Theorem 3]{Ouvres}  is:

\begin{thm}[Serre]
Let $A$ be a principally polarized abelian variety of dimension $n$, defined over a 
number field $K$. Assume that $n=2, 6$ or $n$ is odd and furthermore assume that 
$\mathrm{End}_{\overline{K}}(A)=\mathbb{Z}$. Then there exists a bound $B_{A, K}$ such that, 
for all $\ell>B_{A, K}$, 
\begin{equation*}\mathrm{Im} \overline{\rho}_{A, \ell}=\GSp(A[\ell])\simeq \GSp_{2n}(\mathbb{F}_{\ell}).\end{equation*}
\end{thm}

For arbitrary dimension, the result is not true (see e.g.~\cite{Mumford69} for an example in dimension 4). However, one eventually 
obtains symplectic image by making some extra assumptions. For example, 
there is the following result of C.~Hall (cf.~\cite{Hall}).

\begin{thm}[Hall]\label{thm:Hall}  Let $A$ be a principally polarized abelian 
variety of dimension $n$ defined over a number field $K$, such that $\mathrm{End}_{\overline{K}}(A)=\mathbb{Z}$, and
satisfying the following property:  
\begin{quote} 
{\normalfont (T)} There is a finite extension $L/K$ so that the N\'eron model of $A/L$ over the 
ring of integers of $L$ has a semistable fiber with toric dimension 1.
\end{quote}
Then there is an (explicit) finite constant $B_{A, K}$ such that, for all $\ell\geq B_{A, K}$, 
\begin{equation*}
\mathrm{Im} \overline{\rho}_{A, \ell}=\GSp(A[\ell])\simeq \GSp_{2n}(\mathbb{F}_{\ell}).
\end{equation*} 
\end{thm}

\begin{rem}\label{rem:conditionT}
 In the case when $A=J(C)$ is the Jacobian of a hyperelliptic curve $C$ of genus $n$, 
 say defined by  an equation $Y^2=f(X)$ with $f(X)\in K[X]$ a polynomial of degree $2n+1$, 
 Hall gives a sufficient condition for Condition (T) to be satisfied at a  prime $\mathfrak{p}$ of the ring of integers of $K$; 
 namely, the coefficients of $f(X)$ should have $\mathfrak{p}$-adic valuation greater than or 
 equal to zero and the reduction of $f(X)$ mod $\mathfrak{p}$ (which is well-defined) should 
 have one double zero in a fixed algebraic closure of the residue field, while all the other zeroes are simple.
\end{rem}

Applying the result of Hall with $K=\mathbb{Q}$ yields the following partial answer to the inverse Galois problem:

\begin{cor}\label{cor:invGal}
Let $n\in \mathbb{N}$ be any natural number. Then for \emph{all sufficiently large} primes $\ell$, the group $\GSp_{2n}(\mathbb{F}_{\ell})$ can be realized as a Galois group over $\mathbb{Q}$.
\end{cor}


\begin{rem}
 Several people, including the anonymous referee, pointed us to the following fact: if we consider a family of genus $n$ 
 hyperelliptic curves $C_t$ defined over $\mathbb{Q}(t)$, with big monodromy at $\ell$, then Hilbert's Irreducibility Theorem
 provides us with infinitely many specializations $t=t_0\in \mathbb{Q}$ such that the Jacobian $J_{t_0}$ of the corresponding curve $C_{t_0}$
 satisfies that $\mathrm{Im}\overline{\rho}_{J_{t_0}, \ell}\simeq \mathrm{GSp}_{2n}(\mathbb{F}_{\ell})$. Such families of curves 
 exist for any odd $\ell$ (see e.g.~\cite{Hall2008} or \cite{Zarhin}). In particular, for any $n\in \mathbb{N}$ and any 
 odd $\ell$, the Inverse Galois problem has an affirmative answer for the group $\mathrm{GSp}_{2n}(\mathbb{F}_{\ell})$.
 Although ensuring the existence of the desired curve, this fact does not tell us how to find such a curve explicitly.
\end{rem}

In the case of curves of genus 2,  Le Duff has studied the image of the Galois representations 
attached to the $\ell$-torsion of $J(C)$, when Condition (T) in Theorem \ref{thm:Hall} is satisfied. 
The main result in \cite{LeDuff} is the following:

\begin{thm}[Le Duff]\label{thm:LeDuff}
Let $C$ be a genus $2$ curve defined over $\mathbb{Q}$, with bad reduction of type (II) 
or (IV) according to the notation in \cite{Liu} at a prime $p$. Let $\Phi_p$ be the group 
of connected components of the special fiber of the N\'eron model of $J(C)$ at $p$. For each prime $\ell$ 
and each prime $q$ of good reduction of $C$, let 
$P_{q, \ell}(X)=X^4 + a X^3 + b X^2 + qaX + q^2\in \mathbb{F}_{\ell}[X]$ be the characteristic 
polynomial of the image under $\overline{\rho}_{J(C), \ell}$ of the Frobenius element at $q$ and let
$Q_{q, \ell}(X)= X^2 + aX + b-2q\in \mathbb{F}_{\ell}[X]$, with discriminants $\Delta_P$ and $\Delta_Q$ respectively. 
 
 Then for all primes $\ell$ not dividing $2 pq \vert \Phi_p\vert$ and such that $\Delta_P$ 
 and $\Delta_Q$ are not squares in $\mathbb{F}_{\ell}$, the image of 
 $\overline{\rho}_{J(C), \ell}$ coincides with $\GSp_{4}(\mathbb{F}_{\ell})$.
 
\end{thm}

Using this result, he obtains a realization of $\GSp_4(\mathbb{F}_{\ell})$ as Galois 
group over $\mathbb{Q}$ for all odd primes $\ell$ smaller than 500000. 

\subsection{Explicit surjectivity result}

A key point in Hall's result is the fact that the image under $\overline{\rho}_{A, \ell}$ 
of the inertia subgroup at the place $\mathfrak{p}$ of $L$ which provides the semistable 
fiber with toric dimension 1 is generated by a  nontrivial transvection (whenever $\ell$ 
does not divide $\mathfrak{p}$ nor the cardinality of the group $\Phi_{\mathfrak{p}}$ 
of connected components of the special fiber of the N\'eron model at $\mathfrak{p}$). 
A detailed proof of this fact can be found in Proposition 1.3 of \cite{LeDuff}.

We expand on this point. Given a finite-dimensional vector space $V$ 
over $\mathbb{F}_{\ell}$,  endowed with a symplectic pairing 
$\langle \cdot, \cdot \rangle:V\times V\rightarrow \mathbb{F}_{\ell}$, 
a transvection is an element $T\in \GSp(V, \langle \cdot, \cdot\rangle)$ 
such that there exists a  hyperplane $H\subset V$ satisfying that the 
restriction $T\vert_H$ is the identity on $H$.  We say that it is a nontrivial transvection if $T$ is not the identity\footnote{We adopt the convention that identity is a transvection so that the set of transvections for a given hyperplane $H$ is a group.}. It turns out 
that the subgroups of $\GSp(V, \langle \cdot, \cdot \rangle)$ that contain 
a nontrivial transvection can be classified into three categories as follows (for a proof, see e.g.~\cite[Theorem~1.1]{ArDiWi}):

\begin{thm}\label{thm:classification} 
Let $\ell\geq 5$ be a prime, let $V$ be a finite-dimensional vector space over $\mathbb{F}_{\ell}$, 
endowed with a symplectic pairing 
$\langle \cdot, \cdot \rangle:V\times V\rightarrow \mathbb{F}_{\ell}$ 
and let $G\subset \GSp(V, \langle \cdot, \cdot\rangle)$ be a subgroup 
that contains a  nontrivial transvection. Then one of the following holds:
\begin{enumerate}
 \item $G$ is reducible.
 \item There exists a proper decomposition $V = \bigoplus_{i\in I} V_i$ of $V$ into equidimensional
non-singular symplectic subspaces $V_i$ such that, for each $g \in G$ and each $i \in I$,
there exists some $j \in I$ with $g(V_i) \subseteq V_j$ and such that the resulting action of
$G$ on $I$ is transitive.
 \item $G$ contains $\Sp(V, \langle \cdot, \cdot\rangle)$.
\end{enumerate}
\end{thm}

\begin{rem}
Assume that $V$ is the $\ell$-torsion group of a principally polarized 
abelian variety $A$ defined over $\mathbb{Q}$ and $\langle \cdot, \cdot \rangle$ 
is the symplectic pairing coming from the Weil pairing.  If 
$G=\mathrm{Im}\overline{\rho}_{A, \ell}$ satisfies the third 
condition in Theorem \ref{thm:classification}, then 
$G=\GSp(V, \langle \cdot, \cdot \rangle)$. Indeed, we have the 
following exact sequence
\begin{equation*}
 1 \rightarrow \Sp(V, \langle \cdot, \cdot \rangle)\rightarrow \GSp(V, \langle \cdot, \cdot \rangle) \rightarrow \mathbb{F}_{\ell}^{\times}\rightarrow 1,
\end{equation*}
where the map 
$m:\GSp(A[\ell], \langle \cdot, \cdot \rangle) \rightarrow \mathbb{F}_{\ell}^{\times}$ 
associates to $M$ the scalar $a$ satisfying that, 
for all $u, v\in V$, $\langle Mu, Mv\rangle=a\langle u, v\rangle$. By Equation 
\eqref{eq:multiplier}, the restriction of $m$ to $\mathrm{Im}(\overline{\rho}_{A, \ell})$ 
coincides with the mod $\ell$ cyclotomic character $\chi_{\ell}$.  We can easily conclude the result using that 
$\chi_{\ell}$ is surjective onto $\mathbb{F}_{\ell}^{\times}$. 
\end{rem}

Even in the favourable case when we know that $\mathrm{Im}(\overline{\rho}_{A, \ell})$ 
contains a nontrivial transvection, we still need to distinguish between the three cases in 
Theorem \ref{thm:classification}. In this paper, we will make use of the following consequence of Theorem \ref{thm:classification} (cf.\ Corollary 2.2 of \cite{ArKa}).

\begin{cor}\label{cor:classification} Let $\ell\geq 5$ be a prime, let $V$ be a finite-dimensional vector space over $\mathbb{F}_{\ell}$, 
endowed with a symplectic pairing 
$\langle \cdot, \cdot \rangle:V\times V\rightarrow \mathbb{F}_{\ell}$ and let $G \subset \GSp(V, \langle \cdot, \cdot\rangle)$ 
be a subgroup containing a nontrivial transvection and an element  whose 
characteristic polynomial is irreducible and which has nonzero trace. 
Then $G$ contains $\Sp(V, \langle \cdot, \cdot\rangle)$.
\end{cor}

In order to apply this corollary in our situation, we need some more 
information on the image of $\overline{\rho}_{A, \ell}$. We will obtain 
this by looking at the images of the Frobenius elements $\mathrm{Frob}_q$ 
for primes $q$ of good reduction of $A$.

More generally, let $A$ be an abelian variety defined over a field $K$ 
and assume that $\ell$ is a prime different from the characteristic of $K$. 
Any endomorphism $\alpha$ of $A$ induces an endomorphism of $A[\ell]$, 
in such a way that the characteristic polynomial of $\alpha$ (which is 
a monic polynomial in $\mathbb{Z}[X]$, see e.g.~$\S 3$, Chapter 3 of 
\cite{Lang} for its definition) coincides, after reduction mod $\ell$, 
with the characteristic polynomial of the corresponding endomorphism 
of $A[\ell]$. In the case when $K$ is a finite field (say of cardinality 
$q$), we can consider the Frobenius endomorphism $\phi_{q}\in \mathrm{End}_K(A)$, induced by 
the action of the Frobenius element $\mathrm{Frob}_q\in \mathrm{Gal}(\overline{K}/K)$. 
Then the reduction mod $\ell$ of the characteristic polynomial of $\phi_q$ 
coincides with the characteristic polynomial of $\overline{\rho}_{A, \ell}(\mathrm{Frob}_q)$. 
This will turn out to be particularly useful in the case when $A=J(C)$ 
is the Jacobian of a curve $C$ of genus $n$ defined 
over $K$, since one can determine the characteristic polynomial of
$\overline{\rho}_{J(C), \ell}(\mathrm{Frob}_q)$  by counting the 
$\mathbb{F}_{q^r}$-valued points of $C$, for $r=1, \dots, n$.

As a consequence, we can state the following result, which will be used in the next section.

\begin{thm}\label{thm:explicitsurjectivity}
 Let $A$ be a principally polarized $n$-dimensional abelian variety defined over $\mathbb{Q}$. 
 Assume that there exists a prime $p$ such that the following condition holds:
 
 \begin{quote} {\normalfont ($T_p$)} The special fiber of the N\'eron model of 
 $A$ over $\mathbb{Q}_p$ is semistable with toric dimension $1$.\end{quote} 
 Denote by $\Phi_p$ the group of connected components of the special fiber 
 of the N\'eron model at $p$. 
 Let $q$ be a prime of good reduction of $A$, let $A_q$ be the special 
 fiber of the N\'eron model of $A$ over $\mathbb{Q}_q$ and let 
 $P_q(X)=X^{2n} + aX^{2n-1} + \cdots + q^n\in \mathbb{Z}[X]$ be the characteristic polynomial of the 
 Frobenius endomorphism acting on $A_q$.
 
 Then for all primes $\ell$ which do not divide $6pq\vert \Phi_p\vert a$  
 and are such that the reduction of $P_q(X)$ mod $\ell$ is irreducible 
 in $\mathbb{F}_{\ell}$, the image of $\overline{\rho}_{A, \ell}$ 
 coincides with $\GSp_{2n}(\mathbb{F}_{\ell})$.
\end{thm}

\begin{rem} The condition that $\ell$ does not divide $a$ corresponds to the Frobenius element having non-zero trace modulo $\ell$. Note that the theorem is vacuous when $a=0$. 
\end{rem}

\section{Galois realization of $\GSp_{2n}(\mathbb{F}_{\ell})$ from a 
hyperelliptic curve of genus $n$}\label{sec:3}

Let $C$ be a hyperelliptic curve of genus $n$ over $\Q$, defined by an
equation $Y^2=f(X)$ where $f(X)\in\mathbb{Q}[X]$ is a polynomial of degree~$2n+1$. 
Let $A=J(C)$ be its Jacobian variety. We assume that $A$ satisfies condition
$(T_p)$ for some prime $p$. In this section we present an algorithm, based on
Theorem~\ref{thm:explicitsurjectivity}, which computes a finite set of 
prime numbers $\ell$ for which the Galois representation $\overline{\rho}_{A, 
\ell}$ has image $\GSp_{2n}(\mathbb{F}_{\ell})$. We apply this procedure to an 
example of a genus~3 a curve using a computer algebra system.

\subsection{Strategy}

First, to apply Theorem~\ref{thm:explicitsurjectivity}, we restrict ourselves 
to hyperelliptic curves of genus $n$ whose Jacobian varieties will satisfy Condition 
($T_p$) for some $p$. Namely, we fix a prime number $p$ and then choose  $f(X) 
\in\mathbb{Z}[X]$  monic  of degree $2n+1$ such that both of the following conditions hold:
\begin{enumerate}
 \item The polynomial $f(X)$ only has simple roots over $\overline{\mathbb{Q}}$, so that 
$Y^2=f(X)$ is the equation of an hyperelliptic curve $C$ over $\mathbb{Q}$.
\item All coefficients of $f(X)$ have $p$-adic valuation greater than or equal 
to zero, and the reduction $f(X)\bmod p$ has one double zero in $\overline{\mathbb{F}}_p$, and its other zeroes are simple. This ensures
that $A=J(C)$ satisfies Condition ($T_p$) (see Remark~\ref{rem:conditionT}).
\end{enumerate}

Any prime of good reduction for $C$ is also a prime of good reduction for its 
Jacobian $A$.  Primes of good reduction for the hyperelliptic curve can be 
computed using the discriminant of Weierstrass equations for $C$
(see~\cite{Lockhart}). In our case, it turns out 
that any prime not dividing the discriminant of $f(X)$ is of good 
reduction for $C$, hence for $A$.

We take such a prime number $q$ of good reduction for $A$. Recall that $P_q(X) \in \mathbb{Z}[X]$ is the characteristic polynomial of the Frobenius endomorphism acting on the fiber $A_q$.

Let $\mathcal{S}_q$ denote the set of prime numbers $\ell$ satisfying the
following 
conditions:
\begin{itemize}
 \item[(i)] $\ell$ divides neither  $6pq\vert \Phi_p\vert$ nor the coefficient of $X^{2n-1}$ in $P_q(X)$,
 \item[(ii)] the reduction of $P_q(X)$ modulo $\ell$ is irreducible in 
$\mathbb{F}_{\ell}$.
\end{itemize}
Note that  if the coefficient of $X^{2n-1}$ in $P_q(X)$ is nonzero, 
condition (i)  rules out only finitely many prime numbers $\ell$, whereas if it vanishes, condition (i) rules out all prime numbers $\ell$. By Theorem~\ref{thm:explicitsurjectivity}, for each $\ell \in \mathcal{S}_q$
the 
representation $\overline{\rho}_{A,\ell}$ is surjective with image 
$\GSp_{2n}(\mathbb{F}_{\ell})$. Also, primes in $\mathcal{S}_q$ can be computed
effectively up to a given fixed bound.

Since we want the polynomial $P_q(X)$ (of degree $2n$) to be irreducible modulo $\ell$, its Galois group $G$ over $\Q$ must be a transitive subgroup of $S_{2n}$ with a ${2n}$-cycle. Therefore, by an application of a weaker version of the Chebotarev density theorem due to Frobenius (\cite{SL}, ``Theorem of Frobenius", p.~32), the density of 
$\mathcal{S}_{q}$ is
\[
 \frac{\# \{ \sigma \in G \subset S_{2n} \colon \sigma \text{ is a $2n$-cycle} \}}{\# G}.
\]
This estimate is far from what Theorem~\ref{thm:Hall} provides us, namely that the 
density of $\ell$'s with $\mathrm{Im}(\overline{\rho}_{A, 
\ell})=\GSp_{2n}(\mathbb{F}_{\ell})$ is $1$.

This leads us to discuss the role of the prime $q$. First of all, we can see that 
\[
\bigcup_{q} \mathcal{S}_q = \{\ell  \textrm{ prime} \colon \ell \nmid 6p|\Phi_p| \mbox{ and } \overline{\rho}_{A,\ell} \textrm{ surjective} \},
\]
where the union is taken over all primes $q$ of good reduction for $A$. Note that the inclusion $\subset$ follows directly from Theorem \ref{thm:explicitsurjectivity}. To show the other inclusion $\supset$,  suppose now  that $\ell \nmid 6p|\Phi_p|$ and that the representation at $\ell$ is surjective.  Its image $\GSp_{2n}(\mathbb{F}_{\ell})$  contains an element  with irreducible characteristic 
polynomial and nonzero trace (see for instance Proposition A.2 of \cite{ArKa}). This element defines a conjugacy class $C\subset \GSp_{2n}(\mathbb{F}_{\ell})$ and the Chebotarev density theorem ensures that there exists  $q$ such that $\overline{\rho}_{A,\ell}(\mathrm{Frob}_q)\in C$, hence $\ell \in \mathcal{S}_q$.

%

Moreover, if, for some fixed $\ell$, the events ``$\ell$ belongs to $\mathcal{S}_q$'' are independent as $q$ varies, 
the density of primes $\ell$ for which $\overline{\rho}_{A,\ell}$ is surjective will increase when we take several different primes $q$. A sufficient condition for this density to tend to 1 is that there exists an infinite family of primes $q$ for which the splitting fields of $P_q(X)$ are pairwise linearly disjoint  over $\mathbb{Q}$.

Therefore, it seems reasonable to expect that computing the sets $\mathcal{S}_q$ for several 
values of $q$ increases the density of primes $\ell$ for which we know the surjectivity of 
$\overline{\rho}_{A,\ell}$. This is what we observe numerically in the next example.


\subsection{A numerical example in genus $3$}

We consider the hyperelliptic curve $C$ of genus $n=3$ over $\Q$ defined by
$Y^2=f(X)$, where
\[
 f(X) = X^2(X-1)(X+1)(X-2)(X+2)(X-3)+7(X-28) \in \mathbb{Z}[X].
\]
This is a Weierstrass equation, which is minimal at all primes $\ell$ different from $2$
(see~\cite[Lemma~2.3]{Lockhart}), with discriminant $-2^{12} \cdot
7 \cdot 73\cdot 1069421 \cdot 11735871491$. 
Thus, $C$ has good reduction away from the primes appearing in this
factorization. Clearly, $p=7$ is a prime for which the reduction of $f(X)$ modulo $7$ has one double zero in $\overline{\mathbb{F}}_7$ and otherwise only simple zeroes. Therefore, its Jacobian $J(C)$ satisfies Condition ($T_7$). As we computed with \textsc{Magma}, the order of the component group $\Phi_7$ is $2$. Recall that $P_q(X)$ coincides with the characteristic polynomial of the Frobenius endomorphism of the
reduced curve $C$ modulo $q$ over $\mathbb{F}_q$.

\medskip
Our method provides no 
significant result for $q\in\{3,5\}$ because for $q=3$ the characteristic 
polynomial $P_q(X)$ is not irreducible in $\mathbb{Z}[X]$ and for $q=5$ it has 
zero trace in $\mathbb{Z}$. So in this example, we first take $q=11$.  The curve has $11, 135$ and $1247$ points over 
$\mathbb{F}_{11}$, $\mathbb{F}_{11^2}$ and $\mathbb{F}_{11^3}$, respectively. The 
characteristic polynomial $P_{11}(X)$ is
\[
 P_{11}(X) = X^6-X^5+7X^4-35X^3+77X^2-121X+1331
\]
and it is irreducible over $\Q$. Its Galois group $G$ has order $48$ and is
isomorphic to the wreath product $S_2 \wr S_3$. This group is the direct product of $3$ copies of $S_2$, on which $S_3$ acts by permutation (see ~\cite[Chapter~4]{JamesKerber}): An element of $S_2 \wr S_3$ can be written as $((a_1,a_2,a_3), \sigma)$, where $(a_1,a_2,a_3)$ denotes an element of the direct product $S_2\times S_2 \times S_2$ and $\sigma$ an element of $S_3$. The group law is defined as follows:
$$((a_1,a_2,a_3),\sigma)((a_1',a_2',a_3'),\sigma')=((a_1,a_2,a_3)(a_1',a_2',a_3')^{\sigma},\sigma\sigma'),$$ 
where $(a_1',a_2',a_3')^{\sigma} = (a'_{\sigma(1)},a'_{\sigma(2)},a'_{\sigma(3)})$. 
One can also view the wreath product $S_2\wr S_3$ as the centralizer of $(12)(34)(56)$ in $S_6$, through an embedding $\psi : S_2 \wr S_3 \rightarrow S_6$ whose image is isomorphic to the so-called Weyl group of type $B_3$ (\cite[4.1.18 and 4.1.33]{JamesKerber}). More precisely, under $\psi$, the image of an element $((a_1,a_2,a_3),\sigma)\in S_2\wr S_3$ is the permutation of $S_6$ that acts on $\{1,2,...,6\}$ as follows: it first permutes the elements of the sets $E_1=\{1,2\}$, $E_2=\{3,4\}$ and $E_3=\{5,6\}$ separately, according to $a_1$, $a_2$ and $a_3$ respectively (identifying $E_2,E_3$ with $\{1,2\} $ in an obvious way) and then permutes the pairs $E_1,E_2,E_3$ according to the action of $\sigma$ on the indices. For example, denoting $S_2=\{\id,\tau\}$,  the image under $\psi$ of $((\tau,\id,\id),(123))$ is the $6$-cycle $(135246)$.

Let us now determine the elements of $S_2\wr S_3$ which map to  $6$-cycles in $S_6$ through the embedding $\psi$. For   an  element  in $S_2 \wr S_3$ to be of order $6$, it has to be of the form $((a_1,a_2,a_3),\gamma)$ with  $\gamma$  a $3$-cycle in $S_3$. 
Now, $\psi$ sends an element $((a_1,a_2,a_3),\gamma)$ where either one or three $a_i$'s are $\id$, to a product of two disjoint $3$-cycles in~$S_6$. So the elements of $S_2\wr S_3$ which are $6$-cycles in $S_6$  are among the eight elements $((\id,\id,\tau),\gamma)$, $((\id,\tau,\id),\gamma)$, $((\tau,\id,\id),\gamma)$  and $((\tau,\tau,\tau),\gamma)$ with $\gamma=(123)$ or $\gamma=(132)$. 
Moreover, \cite[Theorem~4.2.8]{JamesKerber} (see also  \cite[Lemma~3.1]{Gramain} or \cite{taylor})  ensures that these $8$ elements are conjugate. Since  $\psi((\tau,\id,\id),(123))=(135246)$ is a $6$-cycle, we deduce that the $8$ elements listed above are exactly the elements of $S_2\wr S_3$ which are $6$-cycles in $S_6$.

To conclude, the Galois group $G$, viewed as a subgroup of $S_6$, contains exactly $8$ elements that are $6$-cycles. Therefore, the density of $\mathcal{S}_{11}$ is $8/48=1/6$.

\medskip
We can compute $P_q(X)$ using
efficient algorithms available in \textsc{Magma} \cite{magma} or \textsc{Sage}
\cite{sage}, which are based on $p$-adic methods. We found that there are $6891$ prime numbers $11
\leq \ell \leq 500000$ that belong to $\mathcal{S}_{11}$. For these $\ell$, 
we know that the image of $\overline{\rho}_{A,\ell}$ is 
$\GSp_6(\mathbb{F}_{\ell})$, so the groups 
$\GSp_{6}(\mathbb{F}_{\ell})$ are realized as Galois groups arising from the 
$\ell$-torsion of the Jacobian of the hyperelliptic curve $C$. For instance, the first
ten elements of $\mathcal{S}_{11}$ are
\[
47, 71, 79, 83, 101, 113, 137, 251, 269,271.
\]
Also, the proportion of prime numbers $11 \leq \ell \leq 500000$ in
$\mathcal{S}_{11}$ is about 0.1659, which is quite in accordance with the
density obtained from the Chebotarev density theorem.

By looking at polynomials $P_q(X)$ for several primes $q$ of good 
reduction, we are able to significantly improve the known proportion of primes $\ell$, up to a given bound, for which the Galois representation is
surjective. Namely, we computed that 
\[
 \{\ell \text{ prime}, 11 \leq \ell \leq 500000 \} \subseteq 
\bigcup_{11\leq
q 
\leq 571} \mathcal{S}_q.
\]
As a consequence, for any prime $11 \leq \ell \leq 500000$, the group 
$\GSp_{6}(\mathbb{F}_{\ell})$ is realized as a Galois group arising from the 
$\ell$-torsion of the Jacobian of the hyperelliptic curve $C$. This is reminiscent of Le Duff's
numerical data for $\GSp_{4}(\mathbb{F}_l)$ (see Theorem~\ref{thm:LeDuff}).

Combining all of the above suggests that the single hyperelliptic curve $C$ 
might provide a positive answer to the inverse Galois problem for
$\GSp_{6}(\mathbb{F}_{\ell})$ for any prime $\ell \geq 11$.

\bibliography{Bibliography}
\bibliographystyle{alpha}

\end{document}